 \documentclass[final,5p,times,twocolumn]{elsarticle}

\usepackage{amssymb}

\makeatletter
\def\ps@pprintTitle{%
 \let\@oddhead\@empty
 \let\@evenhead\@empty
 \def\@oddfoot{\centerline{\thepage}}%
 \let\@evenfoot\@oddfoot}
\makeatother

\usepackage{xcolor}

\usepackage{graphics} 
\usepackage{epsfig} 
\usepackage{amsmath} 
\usepackage{amsthm}
\usepackage{amssymb}  
\usepackage{subfigure}

\setcitestyle{authoryear,open={(},close={)}}

\newtheorem{theorem}{Theorem}
\newtheorem{remark}{Remark}
\newtheorem{example}{Example}
\newtheorem{definition}{Definition}

\begin{document}

\begin{frontmatter}



\title{\textbf{Time-Delay Systems with Delayed Impulses:\\ A Unified Criterion on Asymptotic Stability}\tnoteref{label0}}
\tnotetext[label0]{This work was supported by the Natural Sciences and Engineering Research Council of Canada (NSERC) under the grant RGPIN-2020-03934. The first author was supported by a fellowship from the Pacific Institute for the  Mathematical Sciences (PIMS).}
%
\author{Kexue Zhang\corref{cor1}}\ead{kexue.zhang@ucalgary.ca}
\cortext[cor1]{Corresponding author.}
\author{Elena Braverman\corref{cor2}}\ead{maelena@ucalgary.ca}
\address{Department of Mathematics and Statistics, University of Calgary, Calgary, Alberta T2N 1N4, Canada}

%
%
%

\begin{abstract}
The paper deals with the global asymptotic stability of general nonlinear time-delay systems with delay-dependent impulses through the Lyapunov-Krasovskii method. We derive a unified stability criterion which can be applied to a variety of impulsive systems. The cases when each of the continuous dynamics and the impulsive component is either stabilizing or destabilizing are investigated. Both theoretically and numerically, we demonstrate that the obtained result is more general than those existing in the literature.
\end{abstract}


\end{frontmatter}

\section{Introduction}\label{Sec1}

{I}{mpulsive} systems naturally arises when the dynamics of a physical phenomenon produces discontinuous trajectories. Such discontinuities, which are normally called impulses, naturally occur when the state of a dynamical system changes abruptly over a negligible period of time. Impulsive differential equations serve as ideal mathematical models for impulsive systems (see, e.g., \cite{VL-DDB-PS:1989,NP-AMS:1995,XL-KZ:2019B}), which have widespread applications in multi-agent consensus (\cite{ICM-AM-AG-AMG:2016}), network synchronization (\cite{NM-MBM-JK-JL-AA:2012}), secure communication (\cite{TY-LOC:1997}), disease control and treatment (\cite{FC-VC-PP:2020}), etc.

Time-delay effects are frequently encountered in impulsive systems, for example, impulsive vaccination of epidemic models (\cite{MS-EI:2011}), impulsive control of dynamical networks (\cite{WA-DP-MF:2010}), impulsive consensus in networks with multi-agents (\cite{ZWL-ZHG-XS-GF:2012}), and impulsive predator-prey models (\cite{JD-KSJ:2013}). As one of the most basic properties to dynamical systems, stability has been investigated intensively for impulsive systems with time delay (see a recent review paper \cite{XY-XL-QX-PD:2018} and references therein). Recent years have witnessed increasing interest in the study of impulsive systems with time-delay effects considered in the impulses. Impulses of this type are usually called delay-dependent or delayed impulses (see, e.g, \cite{AK-XL-XS:2009,XL-KZ:2019}), namely, the impulses depend on the system states at some historical moments. Examples can be found in impulsive control systems. The time delay existing in the impulsive controllers is the time inevitably required to sample, process and transmit the system states from the sensors and then to update the actuators.

Recently, remarkable progress has been achieved on stability analysis and control applications of dynamical systems with various types of delayed impulses. For instance, a synchronization problem of nonlinear systems with impulses involving discrete-delays was studied in \cite{AK-XL-XS:2009}, and an upper bound of the impulse intervals (i.e., the intervals lying between two successive impulse times) was required to guarantee synchronization. Stability problems for nonlinear systems with impulses involving various types of time delays have also been investigated in \cite{XL-XZ-SS:2017,XL-SS:2017,XL-JW:2018,XL-SS-JW:2019} and some interesting average dwell-time (ADT) conditions on the impulse intervals were derived. However, the continuous evolution of these impulsive systems does not take into account of time-delay effects. Stability results for impulsive systems with time delays presenting in both the impulse and the continuous parts can be found in \cite{WHC-WXZ:2011,XL-KZ:2019,KZ:2019,XL-KZ:2018,XL-KZ-WCX:2016}. But all of these results provide uniform bounds for the impulse intervals, and no ADT conditions have been reported. Due to the existence of time delay, the study of ADT conditions on stability of such impulsive systems is to a large extent challenging.

The above discussion motivates us to revisit the stability analysis of nonlinear systems with both impulses and time-delay effects. We are particularly interested in nonlinear time-delay systems subject to delay-dependent impulses. A unified asymptotic stability result is obtained for systems with stabilizing continuous dynamics and destabilizing (or stabilizing) impulses, systems with destabilizing continuous evolution and stabilizing impulses, or systems with marginal stable continuous dynamics or marginal stable impulse effects. The unified stability criterion provides the (reverse) ADT conditions on the impulse time sequences, and it is more general than existing results in the sense that our stability guarantee does not require the uniform lower (and/or upper) bound of the impulse intervals. To verify the effectiveness and demonstrate the less conservativeness of our stability criterion, our theoretical result is applied to a scalar system with impulses involving distributed delays, a linear impulsive system with discrete delays, and a nonlinear impulsive control system with time delay. Corresponding numerical simulations are also provided.

We structure the remainder of this paper as follows. First, we introduce the problem description and some preliminaries in Section~\ref{Sec2}. Our main result is then presented in Section~\ref{Sec3} with the proof provided in~\ref{proof}. Detailed discussions of the main result and comparison with the existing results are also conducted in this section. Three examples are presented with numerical simulations in Section~\ref{Sec5}. Finally, Section~\ref{Sec6} summaries our results and discusses some possible directions for future research.

\emph{{Notation.}} 
Denote $\mathbb{N}$ the set of positive integers, $\mathbb{R}^+$ the set of non-negative real numbers, and $\mathbb{R}$ the set of all reals. $\mathbb{R}^n$ and $\mathbb{R}^{n\times n}$ represent the $n$-dimensional and $n\times n$-dimensional real spaces equipped with the Euclidean norm and the induced matrix norm, respectively, both denoted by $\|\cdot\|$. For $A\in \mathbb{R}^{n\times n}$, we denote $A^T$ the transpose of $A$ and $\lambda_{max}(A)$ the largest eigenvalue of $A$. Denote $I$ the identity matrix with appropriate dimensions. We say a function $\alpha:\mathbb{R}^+\rightarrow\mathbb{R}$ belongs to class $\mathcal{K}_{\infty}$ if it is continuous, strictly increasing, unbounded, and satisfies $\alpha(0)=0$. Given constants $a$ and $b$ with $b>a$, let 
\begin{align*}
\mathcal{PC}([a,b],\mathbb{R}^n)=& \Big\{ \phi:[a,b]\rightarrow\mathbb{R}^n  \mathrel{\big|} \phi \textrm{ is piecewise} \textrm{  right-}\cr
& \textrm{continuous} \Big\}\cr
\mathcal{PC}([a,\infty),\mathbb{R}^n) =&  \Big\{\phi:[a,\infty)\rightarrow\mathbb{R}^n \mathrel{\big|} \phi|_{[a,c]}\in \mathcal{PC}([a,c],\mathbb{R}^n) \cr
 & \textrm{ for any } c>a \Big\}
\end{align*}
where function $\phi|_{[a,c]}$ denotes the restriction of $\phi$ to the closed interval $[a,c]$. For a positive $\tau$, the linear space $\mathcal{PC}([-\tau,0],\mathbb{R}^n)$ is equipped with the norm $\|\cdot\|_{\tau}$ defined as $\|\varphi\|_{\tau}:=\sup_{s\in[-\tau,0]}\|\varphi(s)\|$ where $\varphi\in \mathcal{PC}([-\tau,0],\mathbb{R}^n)$. For the sake of simplicity, $\mathcal{PC}_{\tau}$ is used for $\mathcal{PC}([-\tau,0],\mathbb{R}^n)$ in the rest of this paper. 
Given $t,\tau\in\mathbb{R}^+$ and $x\in\mathcal{PC}([-\tau,\infty),\mathbb{R}^n)$, function $x_{t}\in\mathcal{PC}_{\tau}$ is defined as $x_{t}(s)=x(t+s)$ for $s\in[-\tau,0]$, and function $x_{t^-}\in\mathcal{PC}_{\tau}$ is denoted by
\begin{eqnarray}
x_{t^-}(s)=\left\{\begin{array}{ll}
x(t^-), &\textrm{ if } s=0\cr
x(t+s), &\textrm{ if } s\not=0\nonumber
\end{array}\right..
\end{eqnarray}

\section{Preliminaries}\label{Sec2}


Consider the impulsive time-delay system
\begin{eqnarray}\label{sys}
\left\{\begin{array}{ll}
\dot{x}(t)=f(t,x_t), ~t\not=t_k \cr
\Delta x(t)=g(t,x_{t^-}),~t=t_k\cr
x_{t_0}=\phi
\end{array}\right.
\end{eqnarray}
where $x(t)\in\mathbb{R}^n$, initial function $\phi\in\mathcal{PC}_{\tau}$, $\tau>0$ represents the maximum delay involved in system~\eqref{sys}, and $f,g:\mathbb{R}^+\times\mathcal{PC}_{\tau} \rightarrow\mathbb{R}^n$ satisfy $f(t,0)\equiv g(t,0)\equiv 0$ for all $t\in\mathbb{R}^+$ so that system \eqref{sys} admits the trivial solution. The state jump is depicted as $\Delta x(t):=x(t)-x(t^-)$ with $x(t^-)$ representing the left-hand limits of $x$ at time $t$. 
The impulse time sequence $\{t_k\}_{k\in\mathbb{N}}$ is strictly increasing and satisfies $\lim_{k\rightarrow\infty}t_k=\infty$. Throughout this paper, we suppose the state $x$ is right continuous at each impulse time and assume $f$ and $g$ satisfy all the necessary conditions (see \cite{GB-XL:1999} for the fundamental theories of system~\eqref{sys}) so that, for any initial condition $x_{t_0}=\phi\in\mathcal{PC}_{\tau}$, system \eqref{sys} has a unique solution $x(t,t_0,\phi)$ in a maximal time interval $[t_0-\tau,t_0+\Gamma)$, where $0<\Gamma\leq \infty$. We use $N(t,s)$ to indicate the number of impulse moments on the half-closed time interval $(s,t]$ with $t>s\geq t_0$.

\begin{definition}[Global Asymptotic Stability]
System \eqref{sys} is said to be globally asymptotically stable (GAS), if system \eqref{sys} is stable and satisfies $\lim_{t\rightarrow\infty}\|x(t)\|=0$, where $x(t):=x(t,t_0,\phi)$ is the solution of \eqref{sys}.
\end{definition}

Next, we present several concepts regarding Lyapunov candidates. We say a function ${V}:\mathbb{R}^+\times\mathbb{R}^n\rightarrow \mathbb{R}^+$ belongs to class $\mathcal{V}_0$ if, for any $x\in\mathcal{PC}(\mathbb{R}^+,\mathbb{R}^n)$, the composite function $v(t):= {V}(t,x(t))$ belongs to $\mathcal{PC}(\mathbb{R}^+,\mathbb{R}^+)$ and can be discontinuous at some $t^*\in \mathbb{R}^+$ only when $x$ has discontinuities at $t^*$. We say a functional ${V}:\mathbb{R}^+\times \mathcal{PC}_{\tau}\rightarrow\mathbb{R}^+$ belongs to class $\mathcal{V}^*_0$ if, for every $x\in\mathcal{PC}([-\tau,\infty),\mathbb{R}^n)$, {the composite function $v(t):= {V}(t,x_t)$ is continuous on $\mathbb{R}^+$}, and ${V}$ is locally Lipschitz with respect to its second argument. We define the upper right-hand derivative of a Lyapunov functional candidate ${V}:\mathbb{R}^+\times \mathcal{PC}_{\tau}\rightarrow\mathbb{R}^+$ along the trajectory of system \eqref{sys} as
$$\mathrm{D}^+ {V}(t,\psi)=\limsup_{h\rightarrow 0^+}\frac{ {V}\left(t+h,x_{t+h}^{(t,\psi)}\right)- {V}(t,\psi) }{h}$$
where $x^{(t,\psi)}$ is a solution to \eqref{sys} satisfying the initial condition $x_t=\psi$ with $\psi\in\mathcal{PC}_{\tau}$, and $h>0$ is close to zero so that the open interval $(t,t+h)$ contains no impulse times. Our objective is to use the Lyapunov-Krasovskii method to study the asymptotic stability of system \eqref{sys}.

\section{The Unified Criterion}\label{Sec3}
In this section, our main result is introduced followed by detailed discussions of its sufficient conditions.

\begin{theorem}\label{Th}
Suppose there exist functions ${V}_1\in \mathcal{V}_0$, ${V}_2\in \mathcal{V}^*_0$, and class $\mathcal{K}_{\infty}$ functions $\alpha_1, \alpha_2, \alpha_3$, and constants $c,\sigma\in\mathbb{R}$, $\varrho_1,\varrho_2,\mu\geq 0$, $\lambda>0$, and $\kappa>0$, such that, for all $t\in \mathbb{R}^+$ and $\psi\in \mathcal{PC}_{\tau}$,
\begin{itemize}
\item[(i)] $\alpha_1(\|\psi(0)\|)\leq {V}_1(t,\psi(0))\leq \alpha_2(\|\psi(0)\|)$ and $0\leq {V}_2(t,\psi)\leq \alpha_3(\|\psi\|_{\tau})$;

\item[(ii)] ${V}(t,\psi):={V}_1(t,\psi(0))+{V}_2(t,\psi)$ satisfies
\[
\mathrm{D}^+{V}(t,\psi) \leq -c {V}(t,\psi);
\]

\item[(iii)] ${V}_1(t,\psi(0)+g(t,\psi))\leq \varrho_1 {V}_1(t^-,\psi(0)) +\varrho_2 \sup_{s\in[-\tau,0]}\{{V}_1(t^-+s,\psi(s))\};$

\item[(iv)] ${V}_2(t,\psi)\leq \kappa \sup_{s\in[-\tau,0]}\{ {V}_1(t+s,\psi(s))\};$

\item[(v)] for arbitrary $t > s\geq t_0$, the following inequality holds
\begin{equation}\label{inequality}
-\sigma N(t,s) -(c-\lambda)(t-s)\leq \mu
\end{equation}
where the constant $\sigma$ is defined as follows:
\begin{itemize}
\item if $c>0$ and $\varrho_1\geq 1$, then
\begin{equation}\label{d1}
\sigma=-\ln(\varrho_1+\varrho_2 e^{c\tau}) \leq 0;
\end{equation}

\item if $c>0$ and $\varrho_1<1$ with $\varrho_1+[(1-\varrho_1)\kappa +\varrho_2]e^{c\tau}\geq 1$, then
\begin{equation}\label{d2}
\sigma=-\ln(\varrho_1+[(1-\varrho_1)\kappa +\varrho_2]e^{c\tau}) \leq 0;
\end{equation}

\item if $c>0$ and $\varrho_1<1$ with $\varrho_1+[(1-\varrho_1)\kappa +\varrho_2]e^{c\tau}< 1$, then $\sigma$ is a positive constant satisfying the following inequality
\begin{equation}\label{d3}
\varrho_1 e^{\sigma}+ [(1-\varrho_1)\kappa +\varrho_2]e^{c\tau} e^{\sigma N(t_k,t_k-\tau)}\leq 1
\end{equation}
for all $k\in \mathbb{N}$;

\item if $c\leq 0$ and $\varrho_1< 1$ with $\varrho_1+(1-\varrho_1)\kappa +\varrho_2< 1$, then $\sigma> 0$ and satisfies
\begin{equation}\label{d4}
\varrho_1 e^{\sigma}+ [(1-\varrho_1)\kappa +\varrho_2]  e^{\sigma N(t_k,t_k-\tau)}\leq 1
\end{equation}
for all $k\in \mathbb{N}$.
\end{itemize}
\end{itemize}
Then system \eqref{sys} is GAS.
\end{theorem}


In what follows, we give the interpretations of sufficient conditions in Theorem~\ref{Th}, discuss various combinations of the parameters $c$ and $\sigma$, and then compare Theorem~\ref{Th} with the existing results.

\begin{remark}
The Lyapunov-Krasovskii method is applied in Theorem~\ref{Th}, and condition (ii) characterizes the system's continuous evolution. Positive $c$ implies the continuous dynamics is stabilizing, whereas negative $c$ indicates the destabilizing continuous dynamics. The Lyapunov functional candidate $V$ is partitioned into $V_1$ and $V_2$. The impulse effects on the function $V_1$ is outlined in condition (iii), while $V_2\in \mathcal{V}_0^*$ as a composite function is continuous in $t$ which implies the functional part $V_2$ is indifferent to the impulses. Coefficients $\varrho_1$ and $\varrho_2$ correspond to the impulse effects of the non-delayed states and delayed states on $V_1$, respectively. This condition has been extensively employed for stability analysis of nonlinear systems with delayed impulses (see, e.g., \cite{XL-KZ:2019,KZ:2019} for detailed interpretations of $\varrho_1$ and $\varrho_2$). Condition (iv) describes a relationship between $V_1$ and $V_2$ so that the impulse effects on $V_1$ described in condition (iii) can be carried over to the overall Lyapunov candidate $V$. However, we are not able to derive from conditions (iii) and (iv) the following inequality
\begin{equation}\label{d}
V(t,\psi^*)\leq e^{-\sigma} V(t,\psi) ~\textrm{  for } t\in\mathbb{R}^+ \textrm{ and } \psi\in\mathcal{PC}_{\tau}
\end{equation}
where the function $\psi^*\in\mathcal{PC}_{\tau}$ is defined as follows
\begin{align*}
\psi^*(s)= \left\{\begin{array}{ll}
\psi(s)+g(t,\psi), &~\textrm{ if } s=0 \cr
\psi(s), &~\textrm{ if } s\in[-\tau,0),
\end{array}\right.
\end{align*}
due to time-delay effects in the impulses (inequality \eqref{d} is a direct generalization of (4b) in \cite{JPH-DL-ART:2008} for nonlinear systems with delay-free impulses to systems involving delayed impulses). We can observe from the proof of Theorem~\ref{Th} in~\ref{proof} that the constant $\sigma$ defined in condition (v) shares an identical role to the constant $d$ in \cite{JPH-DL-ART:2008}: positive $\sigma$ means the impulses are stabilizing, while negative $\sigma$ corresponds to the destabilizing impulses. For systems with different types of continuous dynamics and impulses, parameter $\sigma$ can be determined according to inequalities~\eqref{d1}-\eqref{d4}, respectively. To balance the continuous evolution and the impulse effects, inequality~\eqref{inequality} in condition (v) provides a unified requirement on identifying feasible impulse times  so that the Lyapunov candidate $V$ converges to zero. It can be seen from the proof of Theorem~\ref{Th} that parameter $\lambda$ plays an essential role in ensuring the exponential convergence of the Lyapunov candidate. Then asymptotic stability of system \eqref{sys} can be naturally concluded from condition (i).
\end{remark}


For different combinations of the signs for $c$ and $\sigma$, inequality \eqref{inequality} leads to interesting requirements on the impulse time sequences which are summarized as follows (please refer to \cite{JPH-DL-ART:2008} for a similar discussion for delay-free systems). To secure the GAS of system \eqref{sys}, the impact of time delays on the selection of the impulse time sequences will also be analyzed.

\begin{itemize}
\item If $c>0$ and $\sigma<0$, then the continuous flow of system~\eqref{sys} is asymptotically stable while the impulses can be destabilizing. We must have $c>\lambda$ so that \eqref{inequality} is satisfied. In this case, we can rewrite \eqref{inequality} as follows
\begin{equation}\label{ADT}
N(t,s)\leq \frac{t-s}{T^*} + N^* ~ \textrm{ for } t > s\geq t_0
\end{equation}
where 
\begin{equation}\label{TN}
T^*= \frac{|\sigma|}{|c-\lambda|} \textrm{ and } N^*=\frac{\mu}{|\sigma|}.
\end{equation}
This condition falls in exactly with the notion of average dwell-time (ADT) initiated in \cite{JPH-ASM:1999} for switching systems. It tells that the destabilizing impulses cannot happen too often. To be more specific, there exists at most one impulse per interval of length $T^*$ on average. From \eqref{d1} and \eqref{d2}, we conclude that enlarging $\tau$ increases both $|\sigma|$ and $T^*$, in other words, setting $\tau$ large amplifies the destabilizing effects from the impulses.

\item If $c>0$ and $\sigma=0$, then the continuous dynamics is asymptotically stable along with potentially stabilizing impulses, that is, the impulses do not pose negative effects on the stability of the continuous dynamics. Inequality \eqref{inequality} tells that the convergence rate $\lambda$ of $V$ cannot be larger than $c$, and this provides no requirements on the impulse time sequences. Increasing $\tau$ makes $\sigma<0$, then our prior discussion applies here.

\item If $c>0$ and $\sigma>0$, then system~\eqref{sys} has both stabilizing continuous evolution and impulses. Intuitively, the overall system is GAS due to the stabilizing effects of both continuous dynamics and the impulses. However, different requirements on the convergence rate $\lambda$ of the Lyapunov candidate will pose different conditions on the impulse time sequence. If $c\geq \lambda$ (that is, the convergence rate $\lambda$ is not bigger than that of the continuous dynamics over each impulse interval), then system~\eqref{sys} with arbitrary impulse times is GAS which can be verified by the fact that inequality~\eqref{inequality} proposes no conditions on the impulse time sequence. Increasing $\tau$ with $\varrho_1+[(1-\varrho_1)\kappa+\varrho_2]e^{c\tau}<1$ reduces the largest possible value of $\sigma>0$ satisfying \eqref{d3}. Similarly, system \eqref{sys} with arbitrary impulse times is GAS since both $c$ and $\sigma$ are positive. On the other hand, if enlarging $\tau$ leads to $\varrho_1+[(1-\varrho_1)\kappa+\varrho_2]e^{c\tau}\geq 1$, then $\sigma\leq 0$ is defined in \eqref{d2} and our discussions for the previous two scenarios apply here. We conclude from this case that increasing $\tau$ in the impulses can destroy their stabilizing effects.

If $c<\lambda$, then the exponential convergence rate of the Lyapunov functional $V$ is larger than the convergence rate of $V$ over each impulse interval, and we can rewrite inequality \eqref{inequality} as
\begin{equation}\label{ReADT}
N(t,s)\geq \frac{t-s}{T^*}-N^* ~\textrm{ for } t > s\geq t_0
\end{equation}
where $T^*$ and $N^*$ are defined in \eqref{TN}. Inequality~\eqref{ReADT} is called a reverse ADT condition which demands that any interval of length $T^*$ has at least one impulse on average. This condition implies that the stabilizing impulses should occur frequently enough so that the convergence rate of $V$ can be bigger than $c$. Furthermore, from \eqref{ReADT} and \eqref{d3}, we have 
\begin{equation}\label{inequalities}
\frac{1}{\sigma} \ln\bigg(\frac{1-\varrho_1 e^{\sigma}}{[(1-\varrho_1)\kappa+\varrho_2]e^{c\tau}}\bigg) \geq  N(t_k,t_k-\tau) \geq \frac{\tau}{T^*}-N^*
\end{equation}
provided $\varrho_1+[(1-\varrho_1)\kappa+\varrho_2]e^{c\tau}<1$. It can be seen from \eqref{inequalities} that $N(t_k,t_k-\tau)$ is bounded from both above and below. The upper bound derived from \eqref{d3} is required to guarantee the impulses maintain their stabilizing effects, while the lower bound assures that the stabilizing impulses indeed accelerate the stabilizing process of the entire system. Setting $\tau$ large in the impulses decreases the upper bound and increases the lower bound, and then shrinks the set of feasible impulse time sequences for GAS of system \eqref{sys}.

\item If $c<0$ and $\sigma>0$, then system~\eqref{sys} has unstable continuous dynamics with stabilizing impulses. From~\eqref{inequality}, we can obtain  reverse ADT condition \eqref{ReADT} with $T^*$ and $N^*$ given in \eqref{TN}. For this scenario, condition \eqref{ReADT} demands that there are no excessively long  impulse intervals in order for the impulses to overcome the negative effects of the continuous flow on GAS of the entire system, so that the Lyapunov functional $V$ converges exponentially with rate $\lambda$. Moreover, we can obtain the following estimation on $N(t_k,t_k-\tau)$ from \eqref{ReADT} and \eqref{d4}:
\begin{equation}\label{inequalities1}
\frac{1}{\sigma} \ln\bigg(\frac{1-\varrho_1 e^{\sigma}}{(1-\varrho_1)\kappa+\varrho_2}\bigg) \geq  N(t_k,t_k-\tau) \geq \frac{\tau}{T^*}-N^*.
\end{equation}
The upper bound of $N(t_k,t_k-\tau)$ introduced in~\eqref{inequalities1} requires the impulses cannot happen too frequently in order to preserve the stabilizing effects. Therefore, the occurrence of the delayed impulses should be carefully determined according to both~\eqref{ReADT} and~\eqref{inequalities1}. See Examples~\ref{example2} and~\ref{example3} for demonstrations. The discussion herein also applies to the following scenario which also requires reverse ADT condition \eqref{ReADT}.

\item If $c=0$ and $\sigma>0$, then system~\eqref{sys} has stabilizing impulses with marginally stable continuous dynamics. From \eqref{inequality}, we can derive the reverse ADT condition \eqref{ReADT} with $c=0$, which requires that the stabilizing impulses occur frequently enough to guarantee the exponential convergence of $V$ with rate $\lambda$.

\item If $c\leq 0$ and $\sigma\leq 0$, then inequality \eqref{inequality} cannot be true since the left-hand side of \eqref{inequality} goes to $\infty$ as $t-s$ approaches $\infty$. Intuitively, the continuous and the impulse parts are both destabilizing that implies the overall system is unstable.

\end{itemize}


\begin{remark}
Though the ADT condition \eqref{inequality} has been well discussed in \cite{JPH-DL-ART:2008} for impulsive systems without time delay, the obtained results certainly are not valid for impulsive time-delay systems. It should be emphasized that time-delay effects are included in both the continuous and the impulse portions of system \eqref{sys}. In this respect, condition \eqref{inequality} has been generalized in Theorem~\ref{Th} to deal with impulsive systems involving time delays. The ADT conditions \eqref{ADT} and \eqref{ReADT} correspond to the concept of the average impulse interval introduced in \cite{JL-DWCH-JC:2010}. The ADT condition \eqref{ADT} requires the average impulse interval has length not smaller than $T^*$, whereas the reverse ADT condition \eqref{ReADT} demands the length of the average impulse interval is not larger than $T^*$. None of these two conditions imposes uniform upper or lower bound of the impulse intervals $[t_k,t_{k+1})$ for $k\in \mathbb{N}$. Therefore, Theorem~\ref{Th} is more general than the results in \cite{WHC-WXZ:2011,XL-KZ:2019,KZ:2019,XL-KZ:2018,XL-KZ-WCX:2016} in this regard (see the examples for detailed illustrations in the following section).
\end{remark}

\section{Illustrative Examples}\label{Sec5}
Three examples are presented to verify our theoretical result and the discussions in the previous section. The first example illustrates Theorem~\ref{Th} with positive $c$ and negative $\sigma$.
\begin{example}\label{example1}
Consider the scalar impulsive system with both discrete and distributed delays from \cite{XL-KZ:2019}:
\begin{subequations}\label{e1sys}
\begin{align}
\label{e1sys-a}\dot{x}(t)&=-\mathrm{sat}\left( x(t))+a \mathrm{sat}(x(t-\tau) \right), ~t\not=t_k\\
\label{e1sys-b}\Delta x(t)&=b \mathrm{sat}\left(\int^t_{t-\tau}x(s)\mathrm{d}s \right),~t=t_k
\end{align}
\end{subequations}
with parameters $a=0.2$, $b=0.25$, $\tau=1$, and the saturation function $\mathrm{sat}(z)=\frac{1}{2}(|z+1|-|z-1|)$ for $z\in\mathbb{R}$.
\end{example}

\begin{figure}[t]
\centering
\includegraphics[width=2.5in]{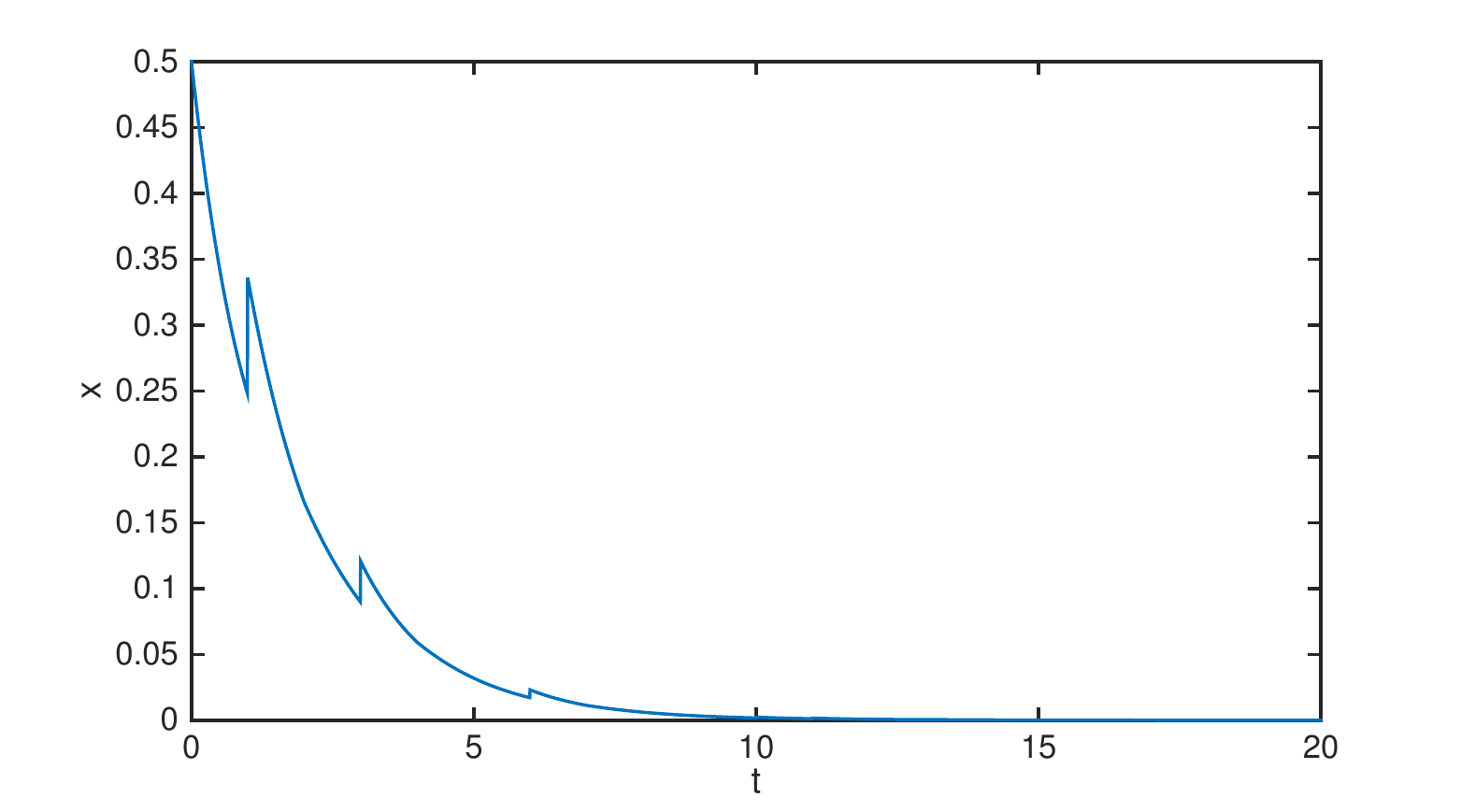}
\caption{Simulation result for system \eqref{e1sys} with $t_0=0$ and initial condition $x_{t_0}(s)=0.5$ for $s\in[-\tau,0]$.}
\label{fig1.1}
\end{figure}

Choose the Lyapunov candidate $V(t,\psi)=V_1(t,\psi(0))+V_2(t,\psi)$ with
\begin{align*}
&V_1(t,\psi(0)) = \left\{\begin{array}{ll}
\psi^2(0), & |\psi(0)|\leq 1 \cr
e^{2(|\psi(0)|-1)}, & |\psi(0)|>1
\end{array}\right.\cr
&V_2(t,\psi) = |a|\int^0_{-\tau}\mathrm{sat}^2(\psi(s))\Big(\epsilon+1+\frac{\epsilon s}{\tau}\Big)\mathrm{d}s
\end{align*}
where $\epsilon>0$. We will check all the conditions of our main result. We first can obtain that condition (i) of Theorem~\ref{Th} holds with
\begin{align*}
\alpha_1(z)=\alpha_2(z)=\left\{\begin{array}{ll}
z^2, & 0\leq z\leq 1 \cr
e^{2(z-1)}, & z>1.
\end{array}\right.
\end{align*}
For $t\not= t_k$, we can see condition (ii) is satisfied, and $c=\min\left\{2-(\epsilon+2)|a|,~\frac{\epsilon}{(\epsilon+1)\tau}\right\}$. This discussion is similar to that in Example 4.1 of \cite{JL-XL-WCX:2011} and thus omitted. To verify condition (iii), we denote $W_1(t)=V_1(t,x(t))$ and $W_2(t)=V_2(t,x_t)$, and then consider two scenarios of $x$ at $t=t_k$.

If $|x(t_k)|\leq 1$, then
\[
W_1(t_k) \leq 2 x^2(t^-_k) + \frac{1}{8} \mathrm{sat}^2 \left( \int^{t_k}_{t_k-\tau}x(s)\mathrm{d}s \right)
\]
and
\begin{align*}
x^2(t^-_k) \leq \left\{\begin{array}{ll}
x^2(t^-_k), & |x(t^-_k)| \leq 1 \cr
e^{2(|x(t^-_k)|-1)}, & |x(t^-_k)|>1,
\end{array}\right.
\end{align*}
that is, $x^2(t^-_k) \leq W_1(t^-_k)$. Identically, we obtain $x^2(t^-_k+s) \leq W_1(t^-_k+s)$ for $s\in [-\tau,0]$, and
\[
\mathrm{sat}^2 \left( \int^{t_k}_{t_k-\tau} x(s)\mathrm{d}s \right) \leq \tau^2 \sup_{s\in[-\tau,0]}\{W_1(t^-_k+s)\}.
\]
Thus, we have
\begin{equation}\label{tk1}
W_1(t_k)\leq 2 W_1(t^-_k) +\frac{\tau^2}{8} \sup_{s\in[-\tau,0]}\{W_1(t^-_k+s)\}.
\end{equation}

If $|x(t_k)|> 1$, then
\begin{equation}\label{tk2}
W_1(t_k)\leq 2e W_1(t^-_k).
\end{equation}
A similar discussion of \eqref{tk2} can be found in \cite{XL-KZ:2019} with more details.

Therefore, we conclude from \eqref{tk1} and \eqref{tk2} that condition (iii) is satisfied with $\varrho_1=2e > 1$ and $\varrho_2=1/8$. Let $\epsilon=4$, then $c=0.8>0$. From Theorem~\ref{Th}, we have $\sigma=-\ln(\varrho_1+\varrho_2 e^{c\tau})\approx -1.7431$. Since $c>0$ and $\sigma<0$, we have that the ADT condition \eqref{ADT} holds with $T^*>|\sigma|/c\approx 2.1789$.

To verify Theorem~\ref{Th} with $c>0$ and $\sigma<0$, we consider system \eqref{e1sys} with the following impulse times:
$t_{4k-3}  = 10k-9$, $t_{4k-2}  = 10k-7$, $t_{4k-1}  = 10k-4$, $t_{4k}   =  10k$, for $k\in\mathbb{N}$. It can be observed that \eqref{ADT} holds with $\mu\geq 4|\sigma|$. Fig. \ref{fig1.1} shows the GAS property of system \eqref{e1sys}. However, Theorem~1 in \cite{XL-KZ:2019} requires $t_k-t_{k-1}\geq 2.1789$ for all $k\in \mathbb{N}$. Since $t_{4k+1}-t_{4k}=1$ and $t_{4k+2}-t_{4k+1}=2$ are both smaller than this lower bound, the result in \cite{XL-KZ:2019} is not applicable to system \eqref{e1sys} with the given impulse times.


In the next example, GAS of a linear impulsive system with discrete delays is studied. Three simulation results are provided, respectively, according to the following combinations of coefficients $c$ and $\sigma$: (\textbf{C1}) $c>0$ and $\sigma<0$;  (\textbf{C2}) $c>0$ and $\sigma>0$; (\textbf{C3}) $c<0$ and $\sigma>0$. We consider $\varrho_1<1$ in (\textbf{C1}) with this example, while the scenario of $\varrho_1>1$ has been illustrated in Example \ref{example1}.

\begin{example}\label{example2} We consider a linear impulsive time-delay system
\begin{subequations}\label{e2sys}
\begin{align}
\label{e2sys-a}\dot{x}(t)&=A x(t)+B x(t-r), ~t\not=t_k\\
\label{e2sys-b}\Delta x(t)&=C x(t) + D x(t-r),~t=t_k
\end{align}
\end{subequations}
where state $x=(x_1,x_2,...,x_n)^T\in \mathbb{R}^n$, matrices $A,B,C,D\in \mathbb{R}^{n\times n}$ and discrete delay $r> 0$.
\end{example}

Consider the Lyapunov candidate $W(t)=W_1(t)+W_2(t)$ with
\[
W_1(t)= x^Tx  \textrm{ and } W_2(t)=\varepsilon \int^t_{t-r} x^T(s)x(s)\mathrm{d}s
\]
where $\varepsilon>0$ is a constant. Then, condition (i) of Theorem~\ref{Th} is true, $\alpha_1(z)=\alpha_2(z)=z^2$ and $\alpha_3(z)=\varepsilon\tau z^2$ for $z\in\mathbb{R}^+$.

For $t\not=t_k$, the continuous dynamics \eqref{e2sys-a} implies 
\[
\dot{W}_1(t) \leq x^T\left(A+A^T + \varepsilon^{-1} B^T B\right) x + \varepsilon x^T(t-r) x(t-r)
\]
and $\dot{W}_2(t)  =  \varepsilon x^Tx- \varepsilon  x^T(t-r) x(t-r)$, then,
\begin{align*}
\dot{W}(t) \leq x^T\left(A+A^T + \varepsilon^{-1} B^T B + \varepsilon I\right) x.
\end{align*}

From the discrete dynamics~\eqref{e2sys-b}, we have 
\begin{align*}
W_1(t_k) &\leq (1+\xi)x^T(t^-_k)(I+C)^T(I+C)x(t^-_k) \cr
         &  ~~~~ + (1+\xi^{-1}) x^T(t_k-r)D^TDx(t_k-r)\cr
         &\leq (1+\xi)\|I+C\|^2 W_1(t^-_k) + (1+\xi^{-1}) \|D\|^2 W_1(t_k-r)
\end{align*}
where constant $\xi>0$ is to be determined. This implies condition (iii) of Theorem~\ref{Th} is satisfied with 
\[
\varrho_1=(1+\xi)\|I+C\|^2 \textrm{ and } \varrho_2=(1+\xi^{-1}) \|D\|^2.
\]
Condition (iv) is also true with $\kappa=\varepsilon\tau$.

If there exists $c\in\mathbb{R}$ such that
\[
A+A^T + \varepsilon^{-1} B^T B + \varepsilon I\leq -c I,
\]
then condition (ii) is satisfied. However, we will choose $\varepsilon=\|B\|$ and replace this condition with 
$c=-\left(\lambda_{max}(A+A^T ) + 2\|B\|\right)$ in the following discussion. 

In Example \ref{example1}, we have studied the case of $\varrho_1\geq 1$. Therefore, we will focus on $\varrho_1<1$ in this example. To determine the sign of $\sigma$, we will need the following estimations
\begin{align}\label{est01}
  & \varrho_1+[(1-\varrho_1)\kappa +\varrho_2]e^{cr}\cr
= & \left(\sqrt{1-\kappa e^{cr}}\|I+C\| +\sqrt{e^{cr}}\|D\|\right)^2 + \kappa e^{cr}
\end{align}
with 
\[
\xi=\frac{\sqrt{e^{cr}}\|D\|}{\sqrt{1-\kappa e^{cr}}\|I+C\|},
\]
and
\begin{align}\label{est02}
  & \varrho_1 e^{\sigma}+[(1-\varrho_1)\kappa +\varrho_2]e^{\sigma N(t_k,t_k-r)}\cr
= & \left(\sqrt{e^{\sigma}-\kappa e^{\sigma N(t_k,t_k-r)}}\|I+C\| +\sqrt{e^{\sigma N(t_k,t_k-r)}}\|D\|\right)^2 + \kappa e^{\sigma N(t_k,t_k-r)}\cr
\end{align}
with 
\[
\xi=\frac{\sqrt{e^{\sigma N(t_k,t_k-r)}}\|D\|}{\sqrt{e^{\sigma}-\kappa e^{\sigma N(t_k,t_k-r)}}\|I+C\|}.
\]

\begin{figure}[t]
\centering
\subfigure[Case \textbf{(C1)}: $c>0$ and $\sigma<0$]{\label{fig2.1}\includegraphics[width=2.5in]{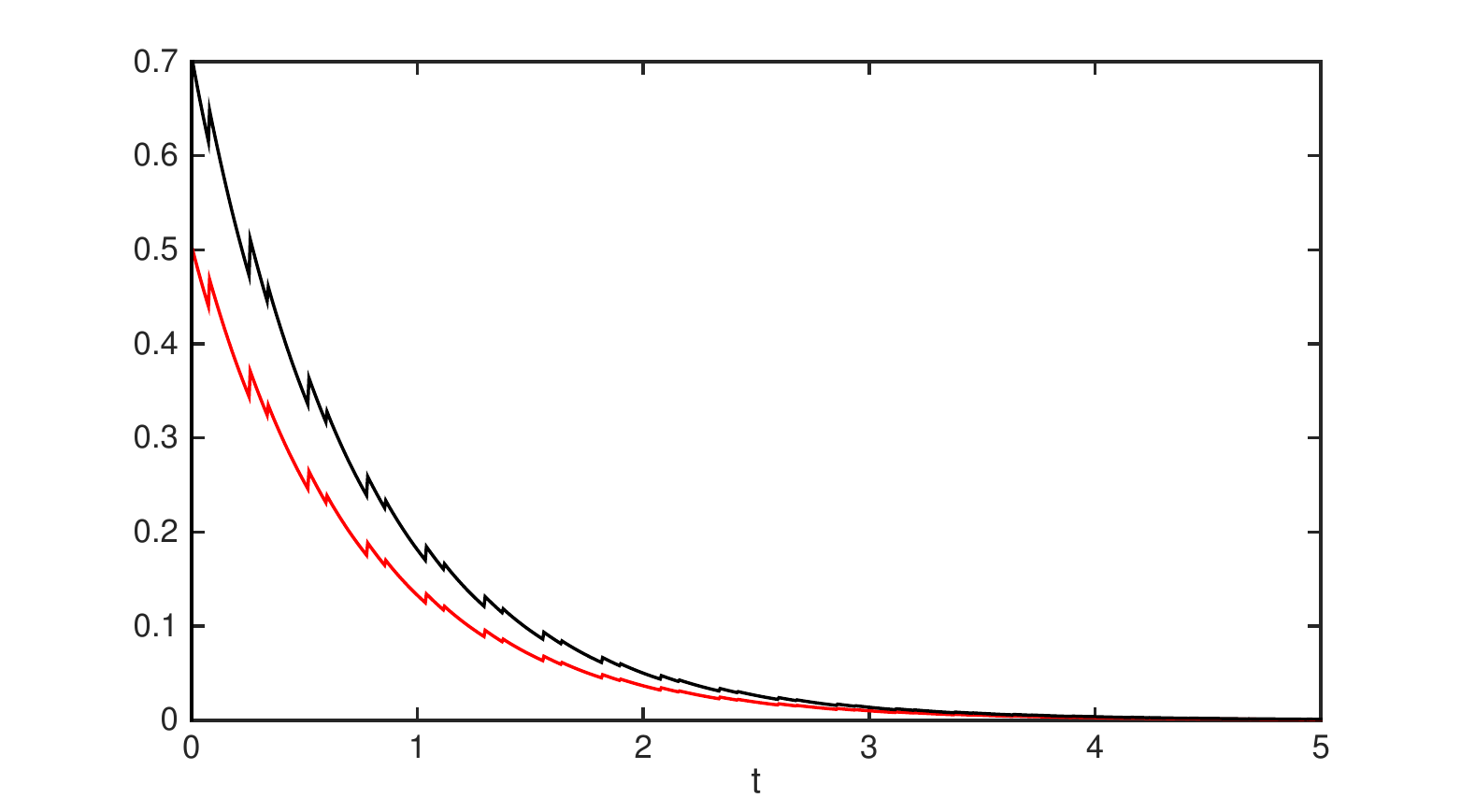}}
\subfigure[Case \textbf{(C2)}: $c>0$ and $\sigma>0$]{\label{fig2.2}\includegraphics[width=2.5in]{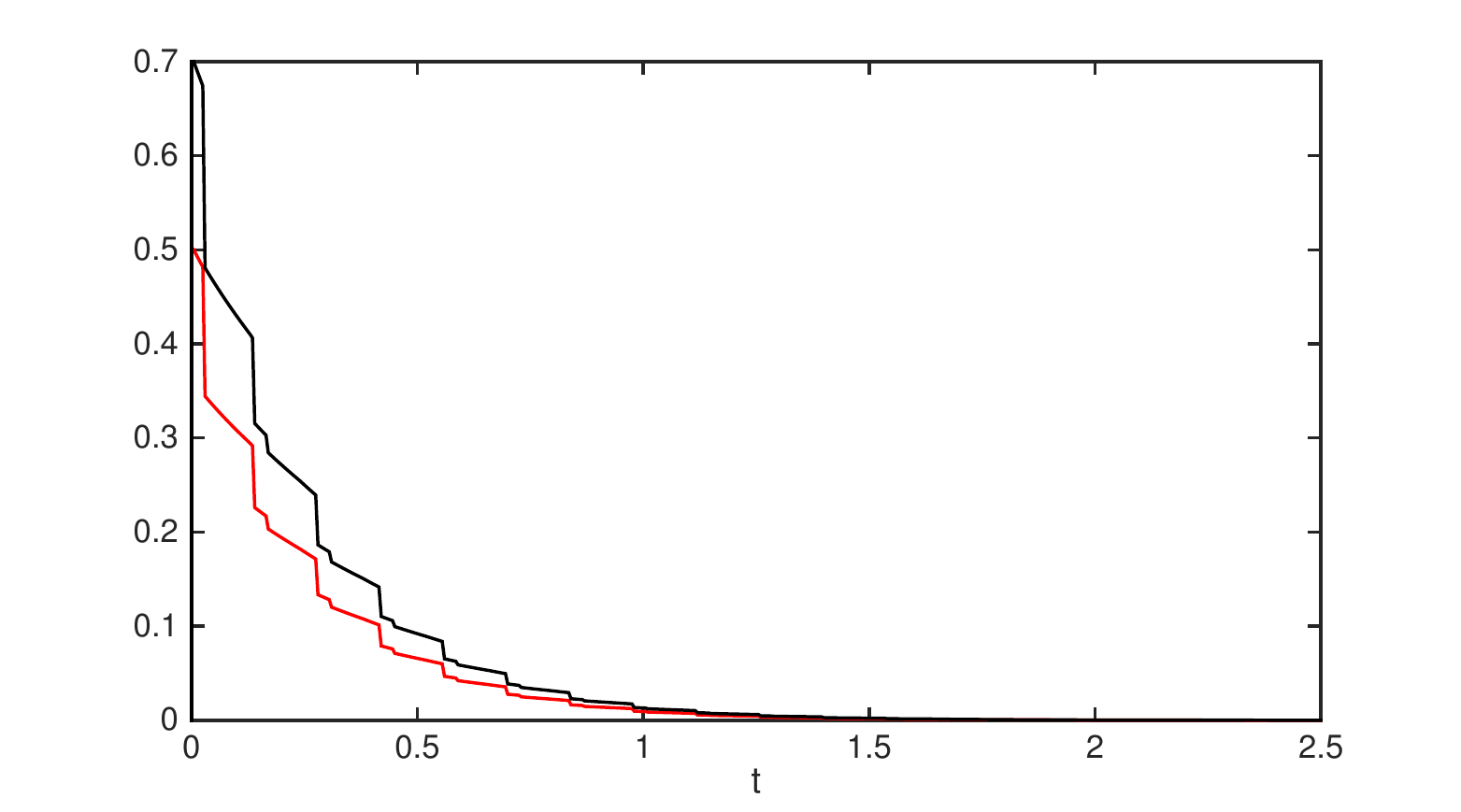}}
\subfigure[Case \textbf{(C3)}: $c<0$ and $\sigma>0$]{\label{fig2.3}\includegraphics[width=2.5in]{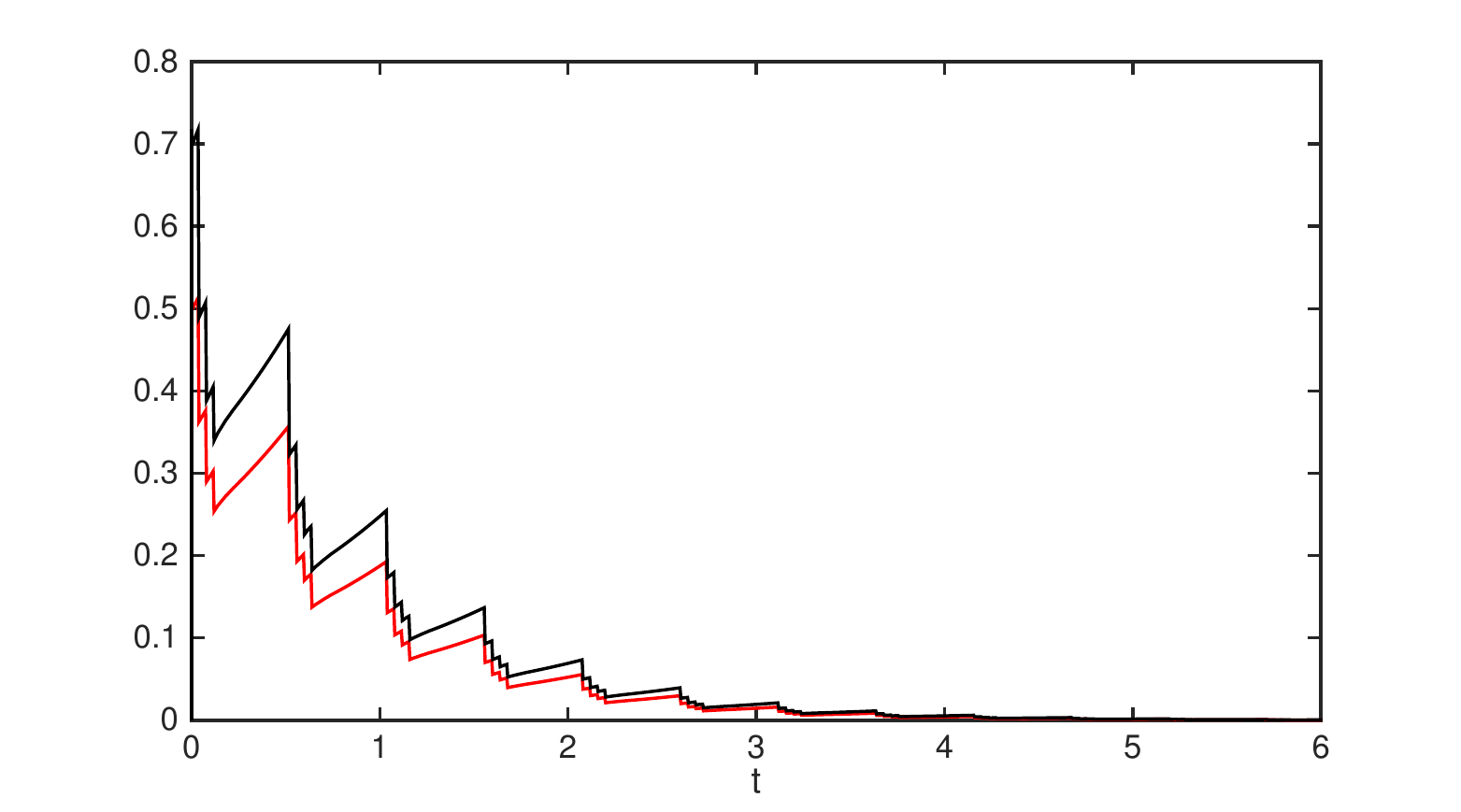}}
\caption{Simulation results for system \eqref{e2sys} with $r=0.1$, $t_0=0$, and initial condition $x_0(s)=[0.5~0.7]^T$ for $s\in [-r,0]$. The red and black curves represent the evolution of states $x_1$ and $x_2$, respectively.}
\label{fig2}
\end{figure}

To verify Theorem~\ref{Th}, we investigate system \eqref{e2sys} with $n=2$ and different coefficients. 

\begin{itemize}
\item[(\textbf{C1})] Consider system \eqref{e2sys} with
\[
A=\begin{bmatrix}
   -1.1834 &  -0.8284 \\
   -0.8284 &  -1.7751
\end{bmatrix}, ~B=\begin{bmatrix}
    0.2500  &  0.1750 \\
    0.1750  &  0.3750
\end{bmatrix},
\]
\[
C=\begin{bmatrix}
    -0.7375  &  0.1750\\
    0.1250  &  -0.6000
\end{bmatrix}, ~D=\begin{bmatrix}
    0.2500  &  0.1750 \\
    0.1750  &  0.3750
\end{bmatrix},
\]
then $\lambda_{max}(A+A^T)\approx-1.1993$, $\|B\|=\|D\|\approx0.4983$, $\|I+C\|\approx0.4972$, and $c\approx 0.2027>0$, $\sigma=-\ln(\varrho_1+[(1-\varrho_1)\kappa +\varrho_2]e^{cr})\approx-0.0262<0$ calculated according to \eqref{est01}. In the simulation, we consider the impulse time sequence: $t_{2k-1}=0.26k-0.18$ and $t_{2k}=0.26k$ for $k\in\mathbb{N}$, then the ADT condition \eqref{ADT} is satisfied with $T^*> |\sigma|/c\approx 0.1293$ and $N^*=\mu/|\sigma|=2$. See Fig. \ref{fig2.1} for an illustration of asymptotic stability of system \eqref{e2sys}. According to Theorem~2 of \cite{XL-KZ:2019}, system \eqref{e2sys} is GAS provided $\inf_{k\in\mathbb{N}}\{t_k-t_{k-1}\}>|\sigma|/c$. Nevertheless, this result cannot be applied to system \eqref{e2sys} with the above given impulse times since $t_{2k}-t_{2k-1}=0.08<|\sigma|/c$.


\item[(\textbf{C2})] Consider system \eqref{e2sys} with matrices $A$, $B$, $D$ given in (\textbf{C1}), and matrix $C$ replaced with the following
\[
C=\begin{bmatrix}
    -0.8950  &  0.0700\\
    0.0500  &  -0.8400
\end{bmatrix},
\]
then $c\approx 0.2027>0$, $\|I+C\|\approx 0.1989$ and $\varrho_1+[(1-\varrho_1)\kappa +\varrho_2]e^{cr}\approx 0.5369 < 1$. For any given impulse time sequence, there exists a uniform lower bound of $t_k-t_{k-1}$ for all $k\in\mathbb{N}$ because of $\lim_{k\rightarrow\infty} t_k=\infty$, which then implies the existence of an upper bound for all $N(t_k,t_k-r)$ with $k\in\mathbb{N}$. Hence, we can find a small enough $\sigma>0$ such that $\varrho_1 e^{\sigma}+[(1-\varrho_1)\kappa +\varrho_2]e^{cr} e^{\sigma N(t_k,t_k-r)}\leq 1$. Since both $c$ and $\sigma$ are positive, we conclude from Theorem~\ref{Th} that system \eqref{e2sys} with an arbitrary impulse time sequence is GAS. The evolution of system \eqref{e2sys} is shown in Fig. \ref{fig2.2} with the following impulse times: $t_{2k-1}=0.14k-0.11$ and $t_{2k}=0.14k$ for $k\in\mathbb{N}$. We can see $\sup_{k\in\mathbb{N}}\{N(t_k,t_k-r)\}=2$, that is, the length of some impulse interval can be smaller than the delay $r$ involved in the impulses.


\item[(\textbf{C3})] Consider system \eqref{e2sys} with matrices
\[
A=\begin{bmatrix}
   0.2 &  0.12 \\
   0.1 &  0.25
\end{bmatrix}, ~D=\begin{bmatrix}
    0.1050  &  0.0700\\
    0.0500  &  0.1600
\end{bmatrix},
\]
$B$ and $C$ given in (\textbf{C1}), then $\lambda_{max}(A+A^T)\approx 0.6756$, $c\approx -1.6722<0$, and $\|D\|\approx 0.1989$. According to \eqref{est02}, inequality $\varrho_1 e^{\sigma}+ [(1-\varrho_1)\kappa +\varrho_2]  e^{\sigma N(t_k,t_k-r)}\leq 1$ holds with $\sigma=0.3786>0$ and $N(t_k,t_k-r)\leq 3$ for all $k\in\mathbb{N}$. Therefore, in the simulation, we consider impulse times $t_{4k-3}=0.52k-0.48$, $t_{4k-2}=0.52k-0.44$, $t_{4k-1}=0.52k-0.4$, $t_{4k}=0.52k$, for $k\in\mathbb{N}$. Then the reverse ADT condition \eqref{ReADT} holds with $T^*< \sigma/|c|\approx 0.2264$ and $N^*=\mu/\sigma=3$, and Theorem~\ref{Th} concludes that system \eqref{e2sys} is GAS (see Fig. \ref{fig2.3} for the numerical simulations). Theorem~3 of \cite{XL-KZ:2019} requires $t_k-t_{k-1}<\ln(\varrho_1+\varrho_2+(1-\varrho_1)\kappa)/c\approx 0.3945$ for all $k\in\mathbb{N}$. Though the upper bound of the impulse intervals obtained from Theorem~3 of \cite{XL-KZ:2019} is bigger than the upper bound of the average impulse interval obtained by our result, Theorem~3 of \cite{XL-KZ:2019} is not applicable to system \eqref{e2sys} with the given impulse time sequence because of $t_{4k}-t_{4k-1}=0.4>0.3945$.
\end{itemize}

\begin{remark}
In the above two examples, the maximum delays in the continuous flows and the impulses are the same. Actually, Theorem~\ref{Th} is applicable to impulsive systems involving different time delays. For example, a general form of system \eqref{e2sys} can be expressed as
\begin{subequations}\label{e3sys}
\begin{align}
\label{e3sys-a}\dot{x}(t)&=A x(t)+B x(t-r_1), ~t\not=t_k\\
\label{e3sys-b}\Delta x(t)&=C x(t) + D x(t-r_2),~t=t_k
\end{align}
\end{subequations}
where $r_1$ and $r_2$ denote, respectively, the time delays in the continuous evolution and the impulses. Define $r:=\max\{r_1,r_2\}$, and it represents the maximum delay in the overall system \eqref{e3sys}, which is in a particular form of system \eqref{sys}. When $r_1=r_2$, system \eqref{e3sys} reduces to system \eqref{e2sys}. The analysis in Example \ref{example2} can be adjusted to study the stability of system \eqref{e3sys} with $r_1\not=r_2$. 
\begin{itemize}
\item If $r_1 > r_2$, then $r=r_1$ and we consider the Lyapunov candidate used in Example \ref{example2}. Similarly to the discussion in Example \ref{example2}, we can derive from the discrete dynamics \eqref{e3sys-b} that
\begin{align*}
W_1(t_k) &\leq \varrho_1 W_1(t^-_k) + \varrho_2 W_1(t_k-r_2)\cr
         &\leq \varrho_1 W_1(t^-_k) + \varrho_2  \sup_{s\in[-r,0]}W_1(t^-_k+s),
\end{align*}
that is, condition (iii) is satisfied. The rest of the discussion is exactly the same as that in Example \ref{example2}.

\item If $0<r_1 < r_2$, then $r=r_2$. Consider the Lyapunov candidate in Example \ref{example2} with $r$ replaced with $r_1$ in $W_2(t)$. The rest discussion is identical to that of Example \ref{example2}.

\item $r_1=0$ implies that the continuous dynamics does not have time delay and $W_2(t)=0$, then $W(t)=W_1(t)$ for all $t\geq t_0$. The analysis of Example \ref{example2} applies to this scenario with $A$ replaced with $A+B$ in the discussion of $\dot{W}_1$. On the other hand, $r_2=0$ means the impulses are free of time delay, and then matrix $D$ can be combined with $C$ in the discussion of $W_1$ at $t=t_k$. 

\end{itemize}

\end{remark}

Last but not least, a nonlinear system with delayed impulses is studied in the following example. Different discrete delays are considered in the continuous dynamics and the impulses of the system.
\begin{example}\label{example3}
We consider the following delayed network control system from~\cite{XL-KZ:2019,JL-XL-WCX:2011}
\begin{subequations}\label{ex3sys}
\begin{align}
\label{ex3sys-a}\dot{x}(t)&=A x(t)+ h(x(t-r)), ~t\not=t_k\\
\label{ex3sys-b}\Delta x(t)&= B x(t-d),~t=t_k
\end{align}
\end{subequations}
where matrices
$A=\begin{bmatrix}
    -18/7 & 9 & 0 \\
    1 & -1 & 1 \\
    0 & -100/7 & 0
\end{bmatrix}$ and $B=-0.5418I$, nonlinear function $h(x)=\mathrm{sat}(x_1)\begin{bmatrix}
    27/7  \\
    0 \\
    0
\end{bmatrix}$ with the saturation function defined in Example~\ref{example1}. In~\eqref{ex3sys-a}, $r=0.02$ represents the time delay in the continuous dynamics. Discrete delay $d=0.01$ in the impulses corresponds to the time required to read the state from the sensors, compute the control input, and update the impulsive actuator. By Theorem~\ref{Th}, we will show that impulse time sequence $\{t_k\}_{k\in\mathbb{N}}$ given as 
\[t_{2k-1}=0.08k-0.03 \textrm{~and~} t_{2k}=0.08k \textrm{~for~} k\in\mathbb{N} \]
guarantees the asymptotic stability of system~\eqref{ex3sys}.
\end{example}

\begin{figure}[t]
\centering
\includegraphics[width=2.5in]{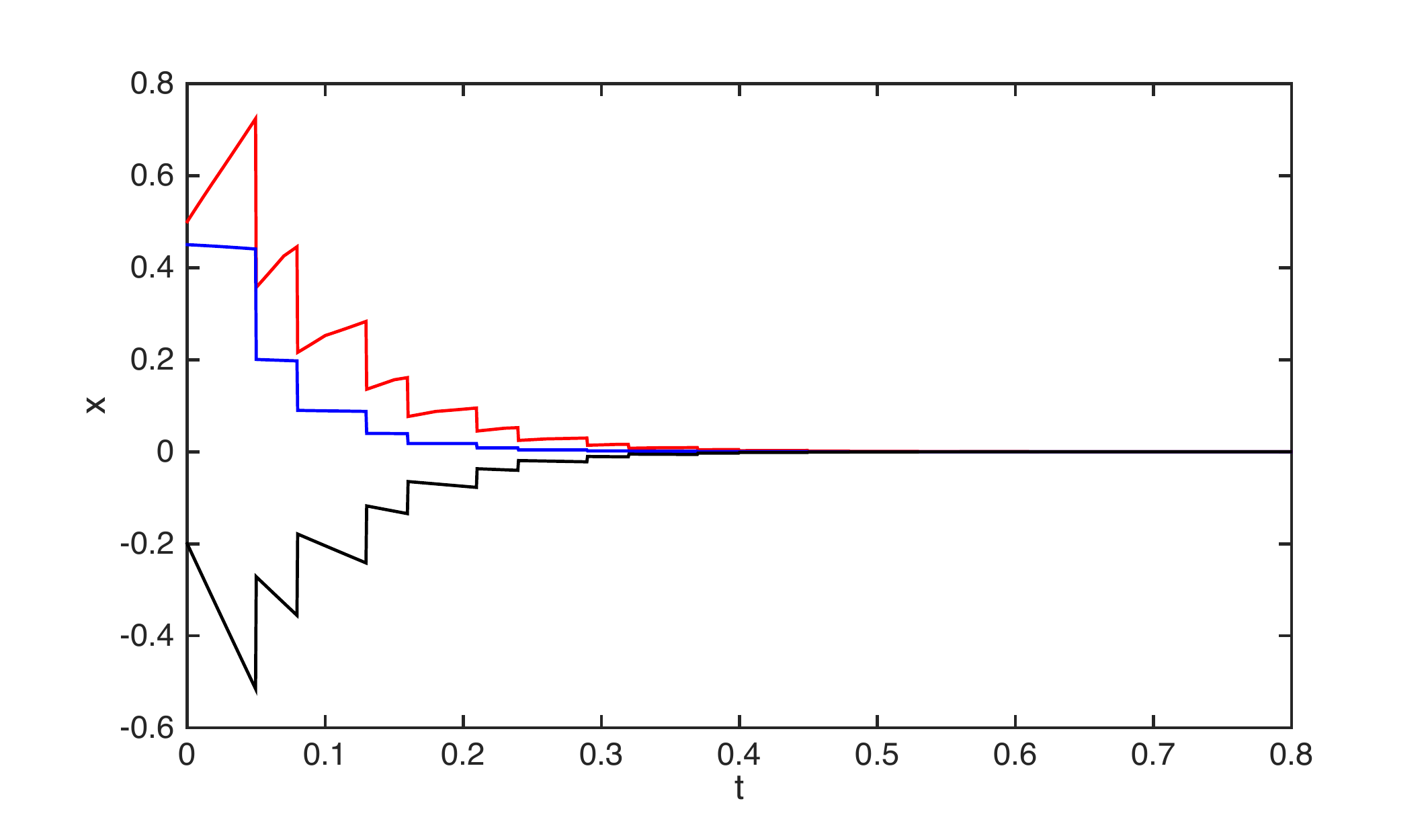}
\caption{State trajectories of system~\eqref{ex3sys} with initial condition $x(s)=[0.5~~0.45~-0.2]^T$ for $s\in [-r,0]$.}
\label{fig3}
\end{figure}

To do so, we use Lyapunov functional candidate $V(t,x_t)=V_1(t,x)+V_2(t,x_t)$ with
\[
V_1(t,x)=x^Tx \textrm{~and~} V_2(t,x_t)= L \int_{t-r}^t x^T(s)x(s)\mathrm{d}s
\]
where $L=27/7$ is the Lipschitz constant of function $h$. Then conditions (i) and (iv) of Theorem~\ref{Th} hold with $\alpha_1(z)=\alpha_2(z)=z^2$, $\alpha_3(z)=rLz^2$ with $z\in\mathbb{R}^+$, and $\kappa=rL$. Similarly to the discussions in~Example 2 from \cite{XL-KZ:2019} and the proof of Theorem 4 in~\cite{XL-KZ-WCX:2016}, we can conclude that conditions (ii) and (iii) of Theorem~\ref{Th} are satisfied with 
\begin{align*}
c &= -\left( \lambda_{max}(A+A^T)+2L \right)<0\cr
\rho_1 &= (1+\xi)\|I+B\|^2 <1 \cr
\rho_2 &= (1+\xi^{-1})\left( d\|B\|(\|A\|+L) +\zeta_k \|B\|^2\right)^2
\end{align*}
where
\[\xi= \frac{d\|B\|(\|A\|+L) +\zeta_k \|B\|^2}{\sqrt{1-\kappa}\|I+B\|},\]
and $\zeta_k=N(t_k-d,t_k)-1$ represents the number of impulses in the open interval $(t_k-d,t_k)$. With the given sequence $\{t_k\}_{k\in\mathbb{N}}$ and impulse delay $d$, we can see that $t_{k+1}-t_k>d$ for all $k\in\mathbb{N}$ which implies that $\zeta_k =0$ for all $k$. Furthermore, $N(t_k-\tau,t_k)=1$ since $t_{k+1}-t_k>\tau=\max\{r,d\}$ for all $k\in\mathbb{N}$. Therefore, inequality~\eqref{d4} holds with $\sigma=0.9619$. We can also observe from the impulse time sequence that the average impulse interval has length $T^*=0.04$, and there exists a small engough constant $\lambda>0$ so that
\[T^*=\frac{|\sigma|}{|c-\lambda|}<\frac{|\sigma|}{|c|}\approx 0.041\]
which then implies inequality~\eqref{d} holds. So far, we have shown all the conditions of Theorem~\ref{Th} are satisfied, and we conclude that system~\eqref{ex3sys} is GAS. State trajectories of system~\eqref{ex3sys} are shown in Fig. \ref{fig3}. Compared with the existing results, Theorem 3 in \cite{XL-KZ:2019} requires $\sup_{k\in\mathbb{N}}\{t_{k+1}-t_k\}<0.041$ whereas our result allows $t_{2k+1}-t_{2k}=0.05>0.041$ for all $k\in\mathbb{N}$.

\section{Conclusions}\label{Sec6}
This paper focused on stability analysis of general nonlinear time-delay systems subject to delayed impulses. We established sufficient conditions on asymptotic stability by using the Lyapunov-Krasovskii functional method. It was shown that the obtained result is more general and applicable to a larger group of impulsive systems than the existing ones. Recently, input-to-state stability (ISS) and integral ISS (iISS) have been studied in \cite{XL-KZ:2019,KZ:2019} for time-delay systems with delayed impulses. Along the lines of this research to investigate ISS and iISS properties of such systems and construct the ADT conditions on the impulse intervals to improve the existing results is a topic for future studies. Another research direction is to generalize our result to hybrid systems with both switchings and delayed impulses. Applications of delayed impulsive control in synchronization and multi-agent consensus of networked systems are also topics for future research.

 \appendix
\section{Proof of Theorem~\ref{Th}} 
\label{proof}
For the sake of notational convenience, we let $W_1(t):={V}_1(t,x(t))$, $W_2(t):={V}_2(t,x_t)$, and then $W(t):=W_1(t)+W_2(t)$.
We use mathematical induction to show
\begin{equation}\label{claim}
W(t)\leq W(t_0) e^{-\sigma N(t,t_0)-c(t-t_0)} \textrm{ for } t\geq t_0.
\end{equation}
When $t\in [t_0,t_1)$, we can derive from condition (ii) that \eqref{claim} holds with $N(t,t_0)=0$.

Next, we suppose \eqref{claim} is true for $t\in[t_0,t_m)$ with some $m\geq 1$ and will show \eqref{claim} still holds on the successive impulse interval $[t_m,t_{m+1})$. When $t=t_m$, we obtain from condition (iii) and the continuity of $W_2$ that
\begin{align}\label{vattm}
W(t_m) & =  W_1(t_m) + W_2(t_m)\cr
       &\leq \varrho_1 W_1(t^-_m) + \varrho_2 \sup_{s\in[-\tau,0]}\{W_1(t^-_m+s)\} + W_2(t^-_m).
\end{align}

To prove \eqref{claim} holds at $t=t_m$, we consider two cases of the constant $c$.

\underline{Case I:} Positive $c$

If $\varrho_1\geq 1$, then we can derive from~\eqref{vattm} and~\eqref{claim} that 
\begin{align*}
W(t_m) & \leq \varrho_1 (W_1(t^-_m)+ W_2(t^-_m) ) + \varrho_2 \sup_{s\in[-\tau,0]}\{W_1(t^-_m+s)\} \cr
       & \leq \varrho_1 W(t^-_m) + \varrho_2 \sup_{s\in[-\tau,0]}\{W(t^-_m+s)\} \cr
       & \leq \varrho_1 W(t_0) e^{-\sigma (N(t_m,t_0)-1)-c(t_m-t_0)} \cr
       &  ~~     + \varrho_2 \sup_{s\in[-\tau,0)}\{ W(t_0) e^{-\sigma N(t_m+s,t_0)-c(t_m+s-t_0)} \}\cr
       & \leq (\varrho_1 + \varrho_2 e^{c\tau}) W(t_0) e^{-\sigma (N(t_m,t_0)-1)-c(t_m-t_0)}\cr
       & \leq W(t_0) e^{-\sigma N(t_m,t_0)-c(t_m-t_0)}
\end{align*}
in which we used the facts that $\sigma\leq 0$, $e^{-\sigma}=\varrho_1 + \varrho_2 e^{c\tau}$ from \eqref{d1}, and $N(t,t_0)\leq N(t_m,t_0)-1$ for $t<t_m$.

If $\varrho_1< 1$ and $\varrho_1+[(1-\varrho_1)\kappa +\varrho_2]e^{c\tau}\geq 1$, then we can derive from~\eqref{vattm},~\eqref{claim}, and condition (iv) that 
\begin{align*}
W(t_m) & \leq \varrho_1 W(t^-_m) + (1-\varrho_1)W_2(t^-_m)+ \varrho_2 \sup_{s\in[-\tau,0]}\{W_1(t^-_m+s)\} \cr
       & \leq \varrho_1 W(t_0) e^{-\sigma (N(t_m,t_0)-1)-c(t_m-t_0)} \cr
       &  ~~ +    [(1-\varrho_1)\kappa+ \varrho_2] \sup_{s\in[-\tau,0)}\{ W(t_0) e^{-\sigma N(t_m+s,t_0)-c(t_m+s-t_0)} \}\cr
       & \leq \left(\varrho_1 + [(1-\varrho_1)\kappa +\varrho_2] e^{c\tau}\right) W(t_0) e^{-\sigma (N(t_m,t_0)-1)-c(t_m-t_0)}\cr
       & \leq W(t_0) e^{-\sigma N(t_m,t_0)-c(t_m-t_0)}.
\end{align*}
In the above inequality, we used the facts $\sigma\leq 0$ and $e^{-\sigma}=\varrho_1 + [(1-\varrho_1)\kappa +\varrho_2] e^{c\tau}$ from \eqref{d2}.

If $\varrho_1< 1$ but $\varrho_1+[(1-\varrho_1)\kappa +\varrho_2]e^{c\tau}< 1$, then we obtain from~\eqref{vattm},~\eqref{claim}, and condition (iv) that 
\begin{align*}
W(t_m) & \leq \varrho_1 W(t^-_m) + (1-\varrho_1)W_2(t^-_m)+ \varrho_2 \sup_{s\in[-\tau,0]}\{W_1(t^-_m+s)\} \cr
       & \leq \varrho_1 W(t_0) e^{-\sigma (N(t_m,t_0)-1)-c(t_m-t_0)} \cr
       &  ~~     +[(1-\varrho_1)\kappa+ \varrho_2] \sup_{s\in[-\tau,0)}\{ W(t_0) e^{-\sigma N(t_m+s,t_0)-c(t_m+s-t_0)} \}\cr
       & \leq \varrho_1 e^{\sigma} W(t_0) e^{-\sigma N(t_m,t_0)-c(t_m-t_0)} \cr
       & ~~ + [(1-\varrho_1)\kappa+ \varrho_2] e^{c\tau}e^{\sigma N(t_m,t_m-\tau)} W(t_0) e^{-\sigma N(t_m,t_0)-c(t_m-t_0)}\cr
       & \leq W(t_0) e^{-\sigma N(t_m,t_0)-c(t_m-t_0)}
\end{align*}
where we used the facts $\sigma>0$ and \eqref{d3}.

\underline{Case II.} Non-positive $c$

If $\varrho_1< 1$ and $\varrho_1+(1-\varrho_1)\kappa+\varrho_2 <1 $, we have
\begin{align*}
W(t_m) & \leq \varrho_1 W(t^-_m) + (1-\varrho_1)W_2(t^-_m)+ \varrho_2 \sup_{s\in[-\tau,0]}\{W_1(t^-_m+s)\} \cr
       & \leq \varrho_1 W(t_0) e^{-\sigma (N(t_m,t_0)-1)-c(t_m-t_0)} \cr
       &  ~~     +[(1-\varrho_1)\kappa+ \varrho_2] \sup_{s\in[-\tau,0)}\{ W(t_0) e^{-\sigma N(t_m+s,t_0)-c(t_m+s-t_0)} \}\cr
       & \leq \varrho_1 W(t_0) e^{-\sigma (N(t_m,t_0)-1)-c(t_m-t_0)} \cr
       & ~~ + [(1-\varrho_1)\kappa+ \varrho_2]  e^{\sigma N(t_m,t_m-\tau)} W(t_0) e^{-\sigma (N(t_m,t_0)-1)-c(t_m-t_0)}\cr
       & \leq W(t_0) e^{-\sigma N(t_m,t_0)-c(t_m-t_0)}.
\end{align*}
Here, we used the facts $\sigma>0$ and \eqref{d4}.

The above discussions conclude that \eqref{claim} holds at $t=t_m$. Then, for $t\in(t_m,t_{m+1})$, we have
\begin{align*}
W(t) & \leq W(t_m) e^{-c(t-t_m)}\cr
     & \leq W(t_0) e^{-\sigma N(t_m,t_0)-c(t_m-t_0)} e^{-c(t-t_m)}\cr
     & =    W(t_0) e^{-\sigma N(t,t_0)-c(t-t_0)}
\end{align*}
where we used $N(t,t_0)=N(t_m,t_0)$ on the impulse interval $(t_m,t_{m+1})$. This completes the proof of the induction. From inequality \eqref{inequality}, we then obtain
\begin{align*}
W(t) & \leq W(t_0) e^{-\sigma N(t,t_0)-c(t-t_0)}\leq e^{\mu}  W(t_0)  e^{-\lambda(t-t_0)}
\end{align*}
for $t\geq t_0$, and applying condition (i) yields
\begin{equation}\label{bdd}
\|x(t)\| \leq \alpha^{-1}_1 \left( e^{\mu}  W(t_0)  e^{-\lambda(t-t_0)} \right).
\end{equation}
This implies that the state $x$ is upper bounded by $\alpha^{-1}_1 \left( e^{\mu}  W(t_0) \right)$ for all $t\geq t_0$. The global existence of solutions to system \eqref{sys} can then be guaranteed (see, \cite{GB-XL:1999}). Therefore, we conclude from \eqref{bdd} that system \eqref{sys} is GAS.






\end{document}